\documentclass[reqno,11pt]{amsart}

\usepackage{xcolor}
\usepackage{epsf}
\usepackage{graphics}
\usepackage{graphicx}
\usepackage{amssymb}
\usepackage{amsmath}
\usepackage{mathcomp}
\usepackage{wasysym}
\date{}

\theoremstyle{plain}
\newtheorem{theorem}{Theorem}
\newtheorem{corollary}{Corollary}

\theoremstyle{definition}
\newtheorem{definition}{Definition}

\theoremstyle{remark}

\def\C{{\mathbb C}}

\def\Z{{\mathbb Z}}

\def\R{{\mathbb R}}

\newcommand{\s}{\sigma}

\newcommand{\scbraid}{{}^{\text{\sl2}^{}}\!B}

\newcommand{\brickgroup}{\langle S_{\tilde\beta} | R_{\tilde\beta}\rangle}
\newcommand{\brickgroupa}{\langle S_{\tilde\alpha} | R_{\tilde\alpha}\rangle}

\title{Secondary braid groups} 

\author{Sebastian Baader \and Michael L\"onne} 

\begin{document}

\begin{abstract}  
We generalize presentations of the fundamental group of discriminant complements
and arrive at a class of presentations associated naturally with words in the free monoid of the alphabet
$\sigma_1,\dots,\sigma_{n-1}$.
Our study addresses invariance properties of these presentations and the presented groups
under various operations on the words.

In particular we prove that the group does only depend on the corresponding element in the positive
braid monoid, and under mild hypotheses only on its conjugacy class in the braid group.
\end{abstract}

\thanks{M.L.\ acknowledges the support by the ERC Grant - 340258 - TADMICAMT }

\maketitle

\section{Introduction}

The aim of this article is to define a group $\scbraid_\beta$ for every positive braid $\beta$.
Our construction will associate a presentation $\langle S_{\tilde\beta}| R_{\tilde \beta} \rangle$ to any word $\tilde \beta$ using letters $\sigma_i$, $i=1,\dots, N-1$ 
and we will study the dependence of the isomorphism type of the associated group on the defining datum. Our definition is motivated by the presentation of the fundamental group of discriminant complements of Brieskorn-Pham polynomials~\cite{L}, while the invariance properties rely on a recent plane graph model for positive braid links~\cite{BLL}. Before even giving detailed definitions, we state our main invariance theorem.

\medskip

\begin{theorem} 
Let $\tilde\alpha$ and $\tilde\beta$ be words in the alphabet
$\sigma_1,\dots,\sigma_{n-1}$, and let
$\alpha,\beta$ be the corresponding positive braids, then
\[
\langle S_{\tilde\alpha}| R_{\tilde \alpha} \rangle
\quad \cong \quad
\langle S_{\tilde\beta}| R_{\tilde \beta} \rangle
\]
if either of the following holds true:
\begin{enumerate}
\item
the positive braids $\alpha,\beta$ are equal.
\item
both $\alpha,\beta$ contain a positive half twist
\[
\Delta=(\sigma_1 \cdots \sigma_{n-1})
(\sigma_1 \cdots \sigma_{n-2}) \cdots (\sigma_1 \cdots \sigma_2) \sigma_1.
\]
and their link closures $\hat\alpha$ and $\hat\beta$ are braid isotopic.
\end{enumerate}
\end{theorem}

We denote the group associated to the
positive braid $\beta$ by $\scbraid(\beta)$ and call it the \emph{secondary braid group}.

\medskip

The group $\scbraid(\beta)$ does not coincide with the fundamental group of the positive braid link $\hat\beta$, in general. As we will shortly see, the groups $\scbraid(\sigma_1^n)$ are isomorphic to the classical braid groups $B_n$, which are not isomorphic to the fundamental group of the corresponding link complements $\widehat{\sigma_1^n}$ except for $n=2,3$. In fact, the group $\scbraid(\beta)$ could be an invariant of positive braid links, a statement that we have not been able to prove so far.

\medskip



A genuine source for the groups $\scbraid(\beta)$ lies in singularity theory:
In case $\tilde\beta$ is a positive word for a braid $\beta$ with closure the torus link of type $T(p,q)$ with $p,q \geq 2$, then the closure is the link of a plane curve singularity $f(x,y)=x^p+y^q \in \C[x,y]$ of Brieskorn-Pham type. The fundamental group of the discriminant complement in its universal unfolding is an invariant of the topological type of $f$ at least for $p<4$, see \cite{L_JST}, and is called the \emph{discriminant knot group} of the singularity $f$.

\begin{corollary} Let $\beta$ be a positive braid on $p$ strands whose closure is braid isotopic to a torus link of type $T(p,q)$ with $p \leq q$. Then the group $\scbraid(\beta)$ is isomorphic to the discriminant knot group of the polynomial $$f(x,y)=x^p+y^q \in \C[x,y].$$
\end{corollary}

Note that 
the conjectural invariance of the discriminant knot group under topological equivalence fits nicely with our theorem.
In fact, each topological type of plane curve singularity determines a unique braid isotopy class of links, such that
any representing positive braid contains the half twist $\Delta$.
\medskip

Related observations corresponding to simple isolated singularities, i.e.\ of type ADE, were already made earlier.
In the context of quivers Grant and Marsh~\cite{GM} have a similar family of groups. One of their main results, Theorem~A, is an invariance theorem for (quiver) braid groups associated to ADE type quivers, up to mutation. Incidentally, in our context, these correspond to positive braid links $\hat{\beta}$ with maximal signature invariant (see \cite{Ba}), i.e.~closures of positive braids whose (symmetrised) Seifert form is positive definite of ADE type on their minimal Seifert surfaces.

\begin{corollary}
Let $\hat{\beta}$ be the closure of a positive braid $\beta$ with maximal signature invariant $\sigma(\hat{\beta})$.
If $\hat\beta$ is prime then the group $\scbraid(\beta)$ is isomorphic to a braid group of Dynkin type ADE.
\end{corollary}

The proofs of Theorem~1 and Corollary~2 base on the fact that positive braid words have braid isotopic closures if and only if they define conjugate braids.
We will use in an essential way the solution of Garside for the conjugacy problem in the braid group and the positive braid monoid to pass through a finite sequence of positive braids using only the positive braid relations and conjugations by permutation braids.

As we will see, the invariance of the presentations $\langle S_{\tilde\beta}| R_{\tilde \beta} \rangle$ under positive Markov moves and under `far commutativity' (i.e. $\sigma_i \sigma_j=\sigma_j \sigma_i$ for $|j-i| \geq 2$) will be an immediate consequence of the definition. 
Therefore, we are left to prove the invariance of the secondary braid groups $\scbraid(\beta)$ under positive braid isotopy, i.e.~under the braid relation and conjugation within positive braids. This is done in Sections~5 and~6, respectively.  In fact we will take care
of elementary conjugation in Section~4, since that allows us to put the location of a braid relation on the top of a braid. 

\medskip

Here we carry on with a definition of the presentations $\langle S_{\tilde\beta}| R_{\tilde \beta} \rangle$ and short proofs of Corollaries~1 and~2. 

\medskip

Let $\tilde\beta$ be a positive word, i.e.~a finite product of letters $\sigma_i$ ($1 \leq i \leq N-1$). The braid diagram on $N$ strands associated to $\tilde\beta$ is the concatenation of the braid diagrams on $N$ strands associated to its letters $\sigma_i$. Brick diagrams were introduced in~\cite{BLL} as a graphical notation for positive braid words. They are obtained by replacing every crossing by a horizontal connection of vertical strands and neglecting loose vertical ends.

\medskip
Brick diagrams ($2$-dimensional) motivated our terminology, since they 
should be seen as the secondary object to the standard line graph ($1$-dimensional) on $n$ vertices. In fact, for the line graph, the edges are the $1$-bricks and correspond to the generators of the standard braid group. The $2$-bricks of the brick diagram of $\tilde\beta$ naturally project to the $1$-bricks of the braid group to which $\beta$ belongs.
The presentations of discriminant knot groups of higher dimensional Brieskorn-Pham singularities, see \cite{L}, naturally fits with a generalization of brick-diagrams to higher dimensions.
\medskip


The main combinatorial object for the definition of secondary braid groups is the \emph{(signed/checkerboarded) linking graph} $\Gamma(\tilde\beta)$, which
is a plane graph together with a sign $+/-$ associated to each of its bounded regions.
It was used in~\cite{BLL} as a model for positive braid links and can be obtained from the plane dual graph of the brick diagram, as follows.
The vertices of $\Gamma(\tilde\beta) \subset \R^2$ are in one-to-one correspondence with bricks, and two vertices are joined by an edge, if and only if their bricks correspond 
to pairs $\s_i\s_i$, $\s_j\s_j$ of one of the following types:

\begin{description}
\item[$\s_i^3$]
$i=j$ and the middle letter $\s_i$ belongs to both bricks,
\item[$\sigma_i\sigma_{i\pm1}\sigma_i\sigma_{i\pm1}$]
$i\pm1=j$ and the different letters alternate in $\beta$.
\end{description}
In that case, we say that the two bricks are \emph{linked}.
If however $\sigma_i\sigma_{i\pm1}\sigma_{i\pm1}\sigma_i$, so the indices are nested, there is no edge between the corresponding vertices and the bricks are not linked, although they share a common boundary segment.

It is readily observed, cf.~\cite{BLL}, that each bounded region of $\Gamma(\tilde\beta)$ is bounded by some positive number of edges of the first type with a common index $i$ throughout and two edges
of the second type with indices $i,i\pm1$. For both of these edges the same sign applies,
which then is the sign of the bounded region.
(The negative sign will be illustrated by a shading in the figures.)

Note that the brick diagram and the linking graph only depend on the isotopy class of the braid diagram 
and do not change under positive Markov moves.

\medskip
In particular,
the use of brick diagrams and plane linking graphs reduces the complexity. In fact, a positive Markov
move does change the braid diagram but neither the brick diagram or the linking graph.
Similarly an isotopy of braid diagrams changes the size of bricks but not the linking graph.
\medskip

Under mild additional assumptions, positive braid links can be classified by a decorated version of linking graphs~\cite{BLL}.
An example of a positive braid word,
$$\tilde\beta=\sigma_1 \sigma_3 \sigma_1 \sigma_2 \sigma_1 \sigma_3 \sigma_1 \sigma_3 \sigma_1 \sigma_2 \sigma_3 \sigma_1 \sigma_3 \sigma_2 
,$$
together with its brick diagram $X(\tilde\beta)$ and linking graph $\Gamma(\tilde\beta)$, is shown in Figure~1.
\begin{figure}[ht]
\scalebox{1.0}{\raisebox{-0pt}{$\vcenter{\hbox{\includegraphics{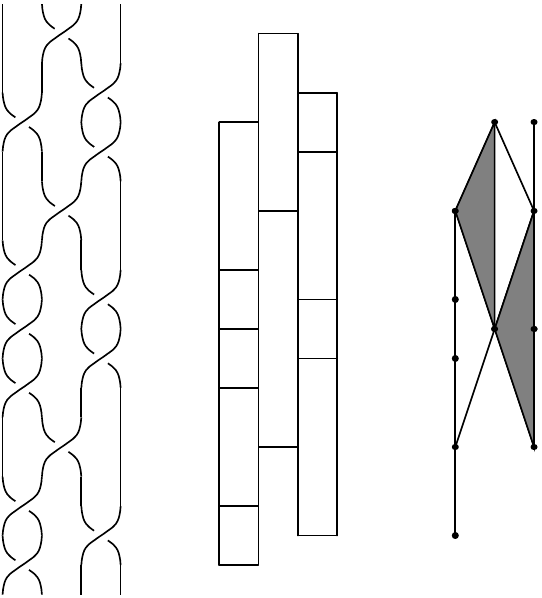}}}$}}
\caption{The braid diagram of $\tilde\beta$, its brick diagram $X(\tilde\beta)$ and linking graph $\Gamma(\tilde\beta)$}
\end{figure}

\begin{definition}
Given a braid positive word $\tilde\beta$ with associated linking graph $\Gamma(\tilde\beta)$.
Let $v_1,\dots, v_k$ be an enumeration of the vertices.
Then the finite presentation $\langle S_{\tilde\beta}| R_{\tilde \beta} \rangle$ is defined in terms of
the \emph{generator set} $S_{\tilde\beta}$ of elements in bijection to the vertices
\[
s_1,s_2,\ldots,s_k
\]
and the \emph{relator set} $R_{\tilde\beta}$ consisting of the following relations:


\begin{description}
\item[braid relation] whenever the corresponding vertices $v_i,v_j$ are joined by an edge.
\[s_i s_j s_i=s_j s_i s_j \]

\item[commutation relation] whenever the corresponding vertices $v_i,v_j$ are not joined by an edge.
\[s_i s_j=s_i s_j \]

\item[cycle relation] whenever the corresponding vertices $v_{i_1},v_{i_2},\ldots,v_{i_n}$
are the vertices in cyclic order around a bounded region of the plane linking graph $\Gamma(\tilde\beta) \subset \R^2$.
\begin{eqnarray*}
&&
s_{i_n} s_{i_{n-1}}\ldots s_{i_2} s_{i_1} s_{i_n} s_{i_{n-1}} \ldots s_{i_4} s_{i_3} \\
&& \qquad = \qquad
s_{i_{n-1}} s_{i_{n-2}}\ldots s_{i_2} s_{i_1} s_{i_n} s_{i_{n-1}} \ldots s_{i_3} s_{i_2}.
\end{eqnarray*}
The cyclic order has to be taken positive or negative according to the sign associated to the region.

\end{description}



\end{definition}






Note that the cycle relation is equivalent to a commutation relation in the presence of the other two kinds of relations:
\begin{eqnarray*}
&&
s_{i_1} \quad s_{i_n} \ldots s_{i_3}s_{i_2}s_{i_3}^{-1} \ldots s_{i_n}^{-1} 
\\ & = &
s_{i_2}^{-1} \ldots s_{i_{n-1}}^{-1}s_{i_n}s_{i_{n-1}} \ldots s_{i_2} \quad s_{i_1}
\\ & = &
s_{i_n} \ldots s_{i_3}s_{i_2}s_{i_3}^{-1} \ldots s_{i_n}^{-1} \quad s_{i_1}
\end{eqnarray*}
Here the second equality is best understood by looking at Figure~5.
A geometric interpretation of the cycle relation is discussed in Section~3 (see also~\cite{Lab,L}).

The easiest example of a secondary braid group is $\scbraid(\sigma_1^2) \cong \Z$. In contrast, the fundamental group of the complement of the closure of the braid $\sigma_1^2 \in B_2$ - the positive Hopf link - is $\Z^2$. Strictly speaking, we did not define the group $\scbraid(\beta)$ for positive braid words $\beta$  without any bricks. These are precisely the positive braids whose closure is a trivial link. For example, the closure of $\sigma_1 \in B_N$ is the trivial link with $N-1$ components. For completeness, we define the secondary braid group of trivial links to be trivial. In the presence of bricks, $\scbraid(\beta)$ admits a surjective homomorphism to $\Z$, by mapping all the generators $s_i$ to $1$. Therefore, secondary braid groups detect trivial links.


By definition, secondary braid groups are quotients of certain Artin groups, as defined by Brieskorn and Saito in~\cite{BS}. The classical braid groups $B_n$ are the secondary braid groups associated with positive braids on two strands:
$$B_n \cong \scbraid(\sigma_1^n).$$
Indeed, $B_n$ is generated by $n-1$ generators $\sigma_1,\sigma_2,\ldots,\sigma_{n-1}$ that commute, except for the pairs with neighbouring indices, which satisfy the braid relation. More generally, if the linking graph $\Gamma(\tilde\beta)$ is a Dynkin diagram of type ADE, then the corresponding group $\scbraid(\beta)$ is a braid group of the same Dynkin type. As a consequence, we obtain a proof of Corollary~2.

\begin{proof}[Proof of Corollary~2]
In \cite{Ba} it is shown that a prime positive braid of braid index at least $4$
either has non-definite Seifert form or can be transformed via braid relations, elementary
conjugation and positive Markov destabilization to a braid word of lower braid index.
If the braid index is at most three and the Seifert form is definite, then
the corresponding braid is shown to be braid isotopic to a braid word with linking graph
of type ADE and the claim follows.
\end{proof}

We proceed with a proof of the first Corollary.

\begin{proof}[Proof of Corollary~1]
Let $\tilde\beta$ be a positive braid word with closure braid isotopic to a torus link of type $T(p,q)$ where $p,q \geq 2$. By Theorem~1, the groups $\scbraid(\beta)$ and $\scbraid(\beta_{(p,q)})$ for $\tilde\beta_{(p,q)}=(\sigma_1 \sigma_2 \ldots \sigma_{p-1})^q$ are isomorphic, since both braids contain a positive half twist and have braid isotopic closures: the torus link $T(p,q)$. 

For the isolated plane curve singularity given by the polynomial $f(x,y)=x^p+y^q \in \C[x,y]$, Hefez and Lazzeri \cite{HL} described the intersection form on the corresponding Milnor lattice with respect to a natural distinguished basis. The
corresponding Coxeter Dynkin diagram was given by Gabrielov \cite{Ga1} and has the same underlying graph as the linking graph $\Gamma(\tilde\beta_{(p,q)})$.
It determines by \cite{L, L_IJM} the presentation of the discriminant knot group of the singularity.
According to Theorem~4.3.\ of \cite{L} the presentation is based on the Coxeter Dynkin diagram in the same way as the presentation of $\scbraid(\beta_{(p,q)})$ is based on the linking graph $\Gamma(\tilde\beta_{(p,q)})$.
Thus the groups $\scbraid(\beta)$ is isomorphic to the discriminant knot group of the singularity $f$.
\end{proof}

Before proceeding with the proof of our theorem, we present a couple of examples that demonstrate the independence of the secondary braid groups $\scbraid(\beta)$ of positive braids $\beta$ and the fundamental groups of their link complements $\pi_1(S^3 \setminus \hat{\beta})$ in Section~2 and discuss the natural action of $\brickgroup$ on the minimal Seifert surface in Section~3.

\section{Examples}

The closures of the two braids $\beta_1,\beta_2 \in B_4$ given by
$$\tilde\beta_1=\sigma_1^4 \sigma_2^3 \sigma_1 \sigma_3 \sigma_2^3 \sigma_3, \; \tilde\beta_2=\sigma_1^4 \sigma_2^2 \sigma_1 \sigma_3 \sigma_2^3 \sigma_3^2$$
are mutant prime knots ($13n1291$ and $13n1320$, in standard notation), hence their fundamental groups are not isomorphic, by a theorem of Gordon and Luecke~\cite{GL}. However, their secondary braid groups are isomorphic, since the linking graphs of $\tilde\beta_1,\tilde\beta_2$ are trees of the same abstract type, see~\cite{BLL}. Here we note that in the absence of cycles in the linking graph $\Gamma(\tilde\beta)$, the isomorphism type of the secondary braid group $\scbraid(\beta)$ does not depend on the planar embedding of $\Gamma(\tilde\beta)$.

The closures of the two braids $\beta_3,\beta_4 \in B_3$ given by
$$\tilde\beta_3=\sigma_1 \sigma_2^2 \sigma_1,\; \tilde\beta_4=\sigma_1 \sigma_2^2 \sigma_1 \sigma_2^2$$
are non-isotopic three-component links with homeomorphic link complements, hence isomorphic fundamental groups. The linking graphs of $\tilde\beta_3,\tilde\beta_4$ are a union of two points and a Dynkin tree of type $D_4$, respectively. The corresponding secondary braid groups are not isomorphic, since they have different abelianisations: $\Z^2$ and $\Z$. 

In the rest of this section, we have a close look at a pair of positive braids with isotopic closures, and construct an isomorphism between their secondary braid groups. We consider the two conjugate positive braids
$\alpha,\beta \in B_3$ given by
$$\tilde\alpha=\sigma_1 \sigma_2 \sigma_1 \sigma_1 \sigma_2 \sigma_1, \; \tilde\beta=\sigma_1 \sigma_1 \sigma_2 \sigma_1 \sigma_1 \sigma_2.$$
The corresponding linking graphs $\Gamma(\tilde\alpha),\Gamma(\tilde\beta)$ are a cycle of length four and the Dynkin diagram $D_4$ (i.e.~a union of three edges sharing one common vertex), respectively, as shown in Figure~2.
\begin{figure}[ht]
\scalebox{1.0}{\raisebox{-0pt}{$\vcenter{\hbox{\includegraphics{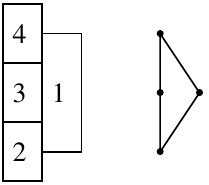}}}$}} \qquad \qquad \scalebox{1.0}{\raisebox{-0pt}{$\vcenter{\hbox{\includegraphics{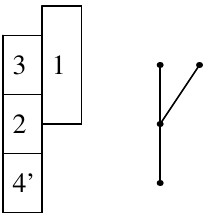}}}$}}
\caption{Linking graphs of braid words $\tilde\alpha$ and $\tilde\beta$}
\end{figure}
Therefore, the cycle relation will come into play in an essential way here. Choosing the enumeration of the four bricks of $\tilde\alpha,\tilde\beta$ as in Figure~2, we obtain the following group presentations:
$$\brickgroupa=\langle s_1,s_2,s_3,s_4; R,R_a \rangle, \; \brickgroup=\langle s_1,s_2,s_3,s'_4; R,R_b \rangle,$$
where $R,R_a,R_b$ are the following sets of relations:
$$R:\; s_1 s_2 s_1=s_2 s_1 s_2,\; s_2 s_3 s_2=s_3 s_2 s_3,\; [s_1,s_3]=1,$$
$$R_a:s_1 s_4 s_1=s_4 s_1 s_4,[s_2,s_4]=1,s_3 s_4 s_3=s_4 s_3 s_4,
s_4 s_3 s_2 s_1 s_4 s_3=s_3 s_2 s_1 s_4 s_3 s_2$$
$$R_b:\; [s_1,s'_4]=1,\; s_2 s'_4 s_2=s'_4 s_2 s'_4,\; [s_3,s'_4]=1.$$
We claim that the map $\varphi: \brickgroupa \to \brickgroup$ taking $s_4$ to $s_3 s_2 s'_4 s_2^{-1} s_3^{-1}$ and leaving the first three generators invariant defines an isomorphism. In fact, we only need to show that $\varphi$ is compatible with the sets of relations $R_a,R_b$, since there is an obvious inverse map $\varphi^{-1}$ taking $s'_4$ to $s_2^{-1} s_3^{-1} s_4 s_3 s_2$. We will do so by substituting
$$\begin{aligned} \varphi^{-1}(s'_4) &= s_2^{-1} s_3^{-1} s_4 s_3 s_2 \\
&= s_2^{-1} s_4 s_3 s_4^{-1} s_2 \\
&= s_4 s_2^{-1} s_3 s_2 s_4^{-1} \\
&= s_4 s_3 s_2 s_3^{-1} s_4^{-1} \\
\end{aligned}$$
for $s'_4$ in the relations of $R_b$ and deriving these from the relations of $R_a$, modulo the relations of $R$.

\begin{enumerate}
\item $[s_3,s'_4]=1$.
$$\begin{aligned} 
s'_4 s_3 &=s_2^{-1} s_3^{-1} s_4 s_3 s_2 s_3=s_2^{-1} s_3^{-1} s_4 s_2 s_3 s_2  \\ 
&=s_2^{-1} s_3^{-1} s_2 s_4 s_3 s_2=s_3 s_2^{-1}  s_3^{-1} s_4 s_3 s_2=s_3 s'_4
\end{aligned}$$
\\
\item $s_2 s'_4 s_2=s'_4 s_2 s'_4$.
$$\begin{aligned} s'_4 s_2 s'_4 &=s_4 s_3 s_2 s_3^{-1} s_4^{-1} s_2 s_4 s_3 s_2 s_3^{-1} s_4^{-1} \\
&=s_4 s_3 s_2 s_3^{-1} s_2 s_3 s_2 s_3^{-1} s_4^{-1}=s_4 s_3 s_2^2 s_4^{-1} \\
&=s_2 s_2^{-1} s_4 s_3 s_4^{-1} s_2^2=s_2 s_2^{-1} s_3^{-1} s_4 s_3 s_2^2=s_2 s'_4 s_2
\end{aligned}$$
\\
\item $[s_1,s'_4]=1$.
$$s'_4 s_1=s_2^{-1} s_3^{-1} s_4 s_3 s_2 s_1$$
$$s_1 s'_4=s_1 s_4 s_3 s_2 s_3^{-1} s_4^{-1} $$
Thus $[s_1,s'_4]=1$ is equivalent to the cycle relation
$$s_4 s_3 s_2 s_1 s_4 s_3=s_3 s_2 s_1 s_4 s_3 s_2.$$
\end{enumerate}

The converse, deriving the relations of $R_a$ from the relations of $R_b$, works analogously by substituting
$$\varphi(s_4)=s_3 s_2 s'_4 s_2^{-1} s_3^{-1}=(s'_4)^{-1} s_2^{-1} s_3 s_2 s'_4$$
for $s_4$ in the relations of $R_b$.

\section{Geometric interpretation of the cycle relation}

The closure $\hat{\beta} \subset S^3$ of a positive braid $\beta \in B_n$ is called a positive braid link. Under the assumption that every generator of $B_n$ appears in $\beta$ at least once, the link $\hat{\beta}$ is non-split and therefore fibred, thanks to a result by Stallings~\cite{St}. This means that the link exterior $S^3 \setminus \hat{\beta}$ admits a fibration over the circle, which is determined by the topological type of one fibre surface~$\Sigma$, together with a gluing map $\varphi: \Sigma \to \Sigma$ called the monodromy map. Thinking of $\R^3$ as $S^3$ minus a point, one can draw an embedded picture of the fibre surface of a positive braid link. Figure~3 shows the fibre surface associated with the braid word $\tilde\alpha=\sigma_1 \sigma_2 \sigma_1 \sigma_1 \sigma_2 \sigma_1$ of the above example. It should be clear from this example that the fibre surface of the closure of a positive braid $\beta$ naturally retracts to the brick diagram $X(\tilde\beta)$ of any representing positive braid word. Therefore the boundaries of the bricks, oriented in the counterclockwise sense, form a basis for the first homology group of the fibre surface.
\begin{figure}[ht]
\scalebox{1.0}{\raisebox{-0pt}{$\vcenter{\hbox{\includegraphics{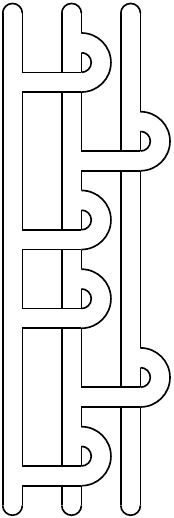}}}$}}
\caption{The fibre surface for the closure of $\alpha=\sigma_1 \sigma_2 \sigma_1 \sigma_1 \sigma_2 \sigma_1$}
\end{figure}

The monodromy map $\varphi: \Sigma \to \Sigma$ admits a decomposition into a product of positive Dehn twist, one for each brick, in a certain order (the order of the bricks starting at the right bottom, going up that column, then proceeding at the bottom of the column to its left, and so on, ending at the top left). This is described in detail in~\cite{Deh}. Going back to the example of the braid $\alpha$, the monodromy map is
$$\varphi=T_{\gamma_4} T_{\gamma_3} T_{\gamma_2} T_{\gamma_1},$$
where $T_{\gamma_i}$ denotes the positive Dehn twist around the curve $\gamma_i$ defined by the brick number~i. We claim that the cycle relation
$$T_{\gamma_4} T_{\gamma_3} T_{\gamma_2} T_{\gamma_1} T_{\gamma_4} T_{\gamma_3}=T_{\gamma_3} T_{\gamma_2} T_{\gamma_1} T_{\gamma_4} T_{\gamma_3} T_{\gamma_2}$$
holds in the mapping class group of the fibre surface $\Sigma$.
In order to see this, we consider the sequence of six curves
$$\gamma_2,T_{\gamma_3}(\gamma_2),T_{\gamma_4}T_{\gamma_3}(\gamma_2),T_{\gamma_1}T_{\gamma_4}T_{\gamma_3}(\gamma_2),T_{\gamma_2}T_{\gamma_1}T_{\gamma_4}T_{\gamma_3}(\gamma_2),T_{\gamma_3}T_{\gamma_2}T_{\gamma_1}T_{\gamma_4}T_{\gamma_3}(\gamma_2)$$
depicted in Figure~4. Setting $h=T_{\gamma_3} T_{\gamma_2} T_{\gamma_1} T_{\gamma_4} T_{\gamma_3}$, we obtain $\gamma_4=h(\gamma_2)$, thus $T_{\gamma_4}=h T_{\gamma_2} h^{-1}$ and $T_{\gamma_4}h=hT_{\gamma_2}$. The latter is nothing but the cycle relation. Needless to say, this argument admits a straightforward generalisation to cycles of arbitrary length.
\begin{figure}[ht]
\scalebox{1.0}{\raisebox{-0pt}{$\vcenter{\hbox{\includegraphics{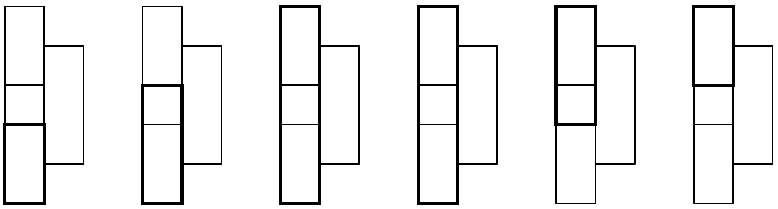}}}$}}
\caption{A sequence of curves relating $\gamma_2$ and $\gamma_4$}
\end{figure}

It is tempting to believe that the group $\scbraid(\beta)$ is isomorphic to the subgroup of the mapping class group of the fibre surface $\Sigma$ of the link $\hat{\beta}$ generated by the Dehn twists around all the curves defined by the bricks. This happens to be true when the linking graph $\Gamma(\beta)$ is a tree of Dynkin type A or D, a fact that was proved by Perron and Vannier~\cite{PV}, based on the work of Birman and Hilden~\cite{BH}. However, this is not true for more complicated tree-like linking graphs. Indeed, Labru\`ere \cite{Lab} proved this for all trees containing a vertex of valency $4$ or two vertices of valency $3$,
while Wajnryb \cite{Waj} showed this for the $E_6$ diagramm, and thus for all trees, that have a vertex of valency $3$
with two neighbours of valency at least $2$.

With this geometric interpretation of the cycle relation at hand, it should not come as a surprise that the cycle relation is equivalent to its various versions with index shifts. More precisely, the two sides of the cycle relation differ by an index shift by one (mod $n$):
$$s_n s_{n-1}\ldots s_2 s_1 s_n s_{n-1} \ldots s_4 s_3=
s_{n-1} s_{n-2}\ldots s_2 s_1 s_n s_{n-1} \ldots s_3 s_2.$$
Any other index shift gives rise to a similar product of $2n-2$ generators $s_i$. The following calculation shows that all these products are equal, provided that one cycle relation holds, as well as the braid and commutation relations. This was already observed by Grant and Marsh (Lemma~2.4.\ in~\cite{GM}). Multiplying both sides of the above relation from the left by $s_1$, and from the right by $s_3^{-1}$, we obtain
$$s_1 s_n s_{n-1}\ldots s_2 s_1 s_n s_{n-1} \ldots s_4=
s_1 s_{n-1} s_{n-2}\ldots s_2 s_1 s_n s_{n-1} \ldots s_3 s_2 s_3^{-1}.$$
The braid relation $s_3 s_2 s_3^{-1}=s_2^{-1} s_3 s_2$ and various commutation relations applied to the right side yield
$$s_1 s_{n-1} s_{n-2}\ldots s_2 s_1 s_2^{-1}s_n s_{n-1} \ldots s_3 s_2.$$
Finally, the braid relation $s_2 s_1 s_2^{-1}=s_1^{-1} s_2 s_1$ and some more commutation relations lead to a new cycle relation, with an index shift by two (mod $n$):
$$s_1 s_n s_{n-1}\ldots s_2 s_1 s_n s_{n-1} \ldots s_4=
s_{n-1} s_{n-2}\ldots s_2 s_1 s_n s_{n-1} \ldots s_3 s_2.$$
All other shifts may be derived by iterating this procedure. From now on, the cycle relation refers to any of its equivalent versions.

\section{Elementary conjugation invariance}

In this section, we prove the invariance of secondary braid groups under elementary conjugation, that is conjugation of the form $\omega \sigma_i \simeq \sigma_i \omega$, where $\omega$ is a positive word and $1 \leq i \leq N-1$. On the level of brick diagrams, this operation amounts to moving one brick in column~$i$ from the top to the bottom. The restriction of the linking graph to that column does not change, but the restriction to the union of the three columns $i-1,i,i+1$ does. Let $n$ be the number of bricks in column~$i$  of the word $\omega \sigma_i$ and let us number them from the bottom to the top. We may assume $n \geq 1$, since otherwise the linking graphs $\Gamma(\omega \sigma_i)$ and $\Gamma(\sigma_i \omega)$ both coincide with $\Gamma(\omega)$. After conjugation, we replace a brick at the top of column~$i$ by a brick at the bottom, which we denote by $n'$ (compare Figure~2 in Section~2). With the example of Section~2 in mind, we define a homomorphism
$$\varphi: {\langle S_{\omega \sigma_i} | R_{\omega \sigma_i}\rangle} \to {\langle S_{\sigma_i \omega} | R_{\sigma_i \omega}\rangle}$$
by mapping the generator $s_n$ to $s_{n-1} s_{n-2} \ldots s_1 s'_n s_1^{-1} \ldots s_{n-2}^{-1} s_{n-1}^{-1}$, and leaving all the other generators invariant. We need to show that this map is compatible with all the relations in $R_{\omega \sigma_i}$ and $R_{\sigma_i \omega}$ involving $s_n$ and $s'_n$, respectively. In fact, this is all we need to show, since there is an obvious inverse map $\varphi^{-1}$ taking $s'_n$ to $s_1^{-1} s_2^{-1} \ldots s_{n-1}^{-1} s_n s_{n-1} \ldots s_2 s_1$.

An important thing to note is the following pair of equalities, which are valid in ${\langle S_{\sigma_i \omega} | R_{\sigma_i \omega}\rangle}$ and ${\langle S_{\omega \sigma_i} | R_{\omega \sigma_i}\rangle}$, respectively. 
$$\begin{aligned} \varphi(s_n) &=s_{n-1} s_{n-2}  \ldots s_1 s'_n s_1^{-1} \ldots s_{n-2}^{-1} s_{n-1}^{-1} \\
&=(s'_n)^{-1} s_1^{-1} \ldots s_{n-2}^{-1} s_{n-1} s_{n-2}  \ldots s_1 s'_n,
\end{aligned}$$
$$\begin{aligned} \varphi^{-1}(s'_n) &=s_1^{-1} s_2^{-1}  \ldots s_{n-1}^{-1} s_n s_{n-1} \ldots s_2 s_1 \\
&=s_n s_{n-1} \ldots s_2 s_1 s_2^{-1} \ldots s_{n-1}^{-1} s_n^{-1}.
\end{aligned}$$

Here we use the fact that the subgroups of ${\langle S_{\sigma_i \omega} | R_{\sigma_i \omega}\rangle}$ and ${\langle S_{\omega \sigma_i} | R_{\omega \sigma_i}\rangle}$ generated by the bricks of column~$i$ are a homomorphic image of an actual braid group $B_{n+1}$ generated by $\s_1,\dots,\s_n$, in which these equalities admit a visual proof, see Figure~5.
\begin{figure}[ht]
\scalebox{1.0}{\raisebox{-0pt}{$\vcenter{\hbox{\includegraphics{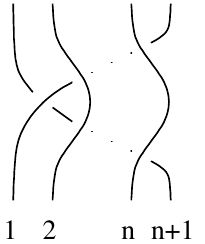}}}$}} \qquad \qquad
\scalebox{1.0}{\raisebox{-0pt}{$\vcenter{\hbox{\includegraphics{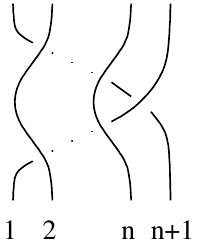}}}$}}
\caption{$\sigma_n \ldots \sigma_2 \sigma_1 \sigma_2^{-1} \ldots \sigma_n^{-1}=\sigma_1^{-1} \ldots \sigma_{n-1}^{-1} \sigma_n \sigma_{n-1} \ldots \sigma_1$}
\end{figure}

Recall that we need to derive all the relations in $R_{\sigma_i \omega}$ involving $s'_n$ from relations in $R_{\omega \sigma_i}$, by substituting $\varphi^{-1}(s'_n)$ for $s'_n$, and vice versa. We start by verifying the braid and commutation relations within the column~$i$.

\begin{enumerate}
\item[(i)] $[s_i,s'_n]=1$, for all $2 \leq i \leq n-1$. This is implied by the fact that the braids of Figure~5 commute with the generators $\sigma_2,\ldots,\sigma_{n-1}$. Here it is essential that $\sigma_{n-1}$ is on the second and third to last strands, rather than the last one.

\item[(ii)] $s_1 s'_n s_1=s'_n s_1 s'_n$. This is equivalent to the fact that the braids of Figure~5 satisfy  the braid relation with the first generator $\sigma_1$. This in turn is equivalent to the equation
$$\sigma_2 \sigma_1 \sigma_2^{-1} \sigma_1 \sigma_2 \sigma_1 \sigma_2^{-1}=\sigma_1 \sigma_2 \sigma_1 \sigma_2^{-1} \sigma_1,$$
since the strands $3$ to $n$ pass straight over all other strands on both sides of $s_1 s'_n s_1$ and $s'_n s_1 s'_n$. Again, a simple picture reveals that both sides of the above equation are equal to $\sigma_2 \sigma_1^2$, so we are done.
\end{enumerate}

Next, we verify all the relations in $R_{\sigma_i \omega}$ involving $s'_n$ and a generator $s_a$ associated with a brick $a$ from another column. Now all relations involve at most two neighbouring columns, except for the commutation relations between bricks that are separated by at least one column, but these obviously transfer from $s_n$ to $s'_n$. Therefore, we may assume that the brick~$a$ is located in column~$i+1$ (the case of column~$i-1$ is symmetric). There are four cases of how the brick~$a$ is linked with the bricks of column~$i$. These cases are shown in Figure~6, before and after the conjugation.
\begin{figure}[ht]
(1) \, \scalebox{0.8}{\raisebox{-0pt}{$\vcenter{\hbox{\includegraphics{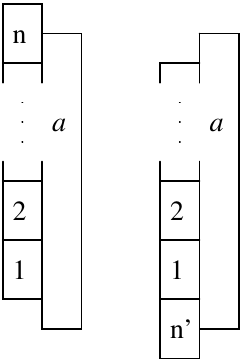}}}$}} \quad
(2) \, \scalebox{0.8}{\raisebox{-0pt}{$\vcenter{\hbox{\includegraphics{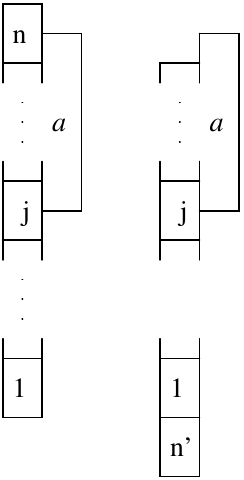}}}$}} \quad
(3) \, \scalebox{0.8}{\raisebox{-0pt}{$\vcenter{\hbox{\includegraphics{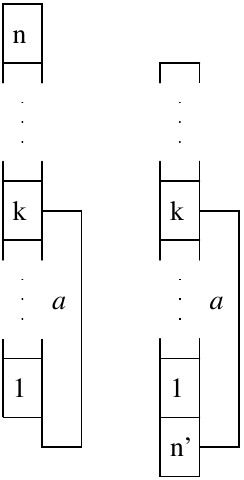}}}$}} \quad
(4) \, \scalebox{0.8}{\raisebox{-0pt}{$\vcenter{\hbox{\includegraphics{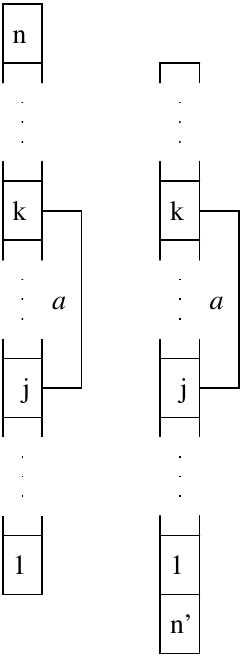}}}$}}
\caption{Four positions for the brick~$a$}
\end{figure}

\begin{enumerate}
\item The brick~$a$ is linked with brick~$n$, and no other brick of column~$i$. In this case, the subgroups of ${\langle S_{\omega \sigma_i} | R_{\omega \sigma_i}\rangle}$ and ${\langle S_{\sigma_i \omega} | R_{\sigma_i \omega}\rangle}$ generated by $s_1,\ldots,s_n,s_a$ and $s_1,\ldots,s'_n,s_a$ are a homomorphic image of the braid group $B_{n+2}$. We need to derive the braid relation $s_a s'_n s_a=s'_n s_a s'_n$ from the braid relation $s_a s_n s_a=s_n s_a s_n$:
$$\begin{aligned}  &\quad\: s'_n s_a s'_n  \\
\qquad&=s_1^{-1} s_2^{-1}  \ldots s_{n-1}^{-1} s_n s_{n-1} \ldots s_2 s_1 s_a s_1^{-1} s_2^{-1}  \ldots s_{n-1}^{-1} s_n s_{n-1} \ldots s_2 s_1 \\
&= s_1^{-1} s_2^{-1}  \ldots s_{n-1}^{-1} s_n s_a  s_n s_{n-1} \ldots s_2 s_1 \\
&= s_1^{-1} s_2^{-1}  \ldots s_{n-1}^{-1} s_a s_n  s_a s_{n-1} \ldots s_2 s_1 \\
&= s_a s_1^{-1} s_2^{-1}  \ldots s_{n-1}^{-1} s_n  s_{n-1} \ldots s_2 s_1 s_a=s_a s'_n s_a. \\
\end{aligned}$$

\item The brick~$a$ is linked with brick~$n$ and one other brick $j$ from column~$i$. We need to derive the commutation relation $[s_a,s'_n]=1$ from the cycle relation for the cycle spanned by the bricks $j,j+1,\ldots,n,a$:
$$s_n \ldots s_j s_a s_n \ldots s_{j+1}=s_{n-1} \ldots s_j s_a s_n \ldots s_{j+1} s_j.$$
For this purpose, we compute
$$A=s_{n-1} \ldots s_1 s'_n s_a  s_1^{-1}  \ldots s_{j-1}^{-1}=s_n s_{n-1} \ldots s_j s_a,$$
$$\begin{aligned}  B &= s_{n-1} \ldots s_1 s_a s'_n  s_1^{-1}  \ldots s_{j-1}^{-1} \\
&=  s_{n-1} \ldots s_j s_a s_j^{-1}  \ldots s_{n-1}^{-1} s_n  s_{n-1} \ldots s_j \\
&=  s_{n-1} \ldots s_j s_a s_n \ldots s_{j+1} s_j s_{j+1}^{-1} \ldots s_n^{-1} \\
\end{aligned}$$ \\
The commutation relation $[s_a,s'_n]=1$ is equivalent to $A=B$, which is equivalent to the cycle relation above.

\item The brick~$a$ is linked with one brick~$k \neq n$, and no other brick from column~$i$. This case is related to case~(2) by a straightforward symmetry; the cycle relation
$$s_k \ldots s_1 s'_n s_a s_k \ldots s_1=s_{k-1} \ldots s_1 s'_n s_a s_k \ldots s_1 s'_n$$
is equivalent to the commutation relation $[s_a,s_n]=1$.

\item The brick~$a$ is linked with two bricks~$j<k<n$ from column~$i$. Here again, we need to derive the commutation relation $[s_a,s'_n]=1$ from the cycle relation for the cycle spanned by the bricks $j,j+1,\ldots,k,a$. We need two steps in order to see this. First, we compute
$$\begin{aligned}  A &= s_{j-1} \ldots s_1 s_a s'_n  s_1 \ldots s_{j-1}=s_a s_j^{-1} \ldots s_{n-1}^{-1} s_n s_{n-1} \ldots s_j \\
&= s_a s_n \ldots s_{j+1} s_j s_{j+1}^{-1} \ldots s_n, \\
\end{aligned}$$
$$\begin{aligned}  B &= s_{j-1} \ldots s_1 s'_n s_a  s_1 \ldots s_{j-1}=s_j^{-1} \ldots s_{n-1}^{-1} s_n s_{n-1} \ldots s_j s_a\\
&= s_n \ldots s_{j+1} s_j s_{j+1}^{-1} \ldots s_n s_a. \\
\end{aligned}$$
Then we compute
$$C=s_{k+1}^{-1} \ldots s_n^{-1} A s_n \ldots s_{k+1}=s_a s_k \ldots s_{j+1} s_j s_{j+1}^{-1} \ldots s_k^{-1},$$
$$\begin{aligned}  D &= s_{k+1}^{-1} \ldots s_n^{-1} B s_n \ldots s_{k+1}=s_k \ldots s_{j+1} s_j s_{j+1}^{-1} \ldots s_k^{-1} s_a \; \\
&= s_j^{-1} \ldots s_{k-1}^{-1} s_k s_{k-1} \ldots s_j s_a. \\
\end{aligned}$$
The commutation relation $[s_a,s'_n]=1$ is equivalent to $A=B$, thus to $C=D$, which is equivalent to the cycle relation
$$s_k \ldots s_j s_a s_k \ldots s_{j+1}=s_{k-1} \ldots s_j s_a s_k \ldots s_j.$$
\end{enumerate}

$$\varphi: {\langle S_{\omega \sigma_i} | R_{\omega \sigma_i}\rangle} \to {\langle S_{\sigma_i \omega} | R_{\sigma_i \omega}\rangle}$$
$$\brickgroup\qquad \brickgroupa \quad R_{\tilde\alpha} \qquad R_{\tilde\beta}$$

\section{Braid relation invariance}

In this section, we verify the invariance of the secondary braid groups under the braid relation. Thanks to the invariance under elementary conjugation established above, it is enough to consider a pair of positive words of the form $\tilde\alpha=\omega \sigma_i \sigma_{i+1} \sigma_i,\tilde\beta=\omega \sigma_{i+1} \sigma_i \sigma_{i+1}$, where $\omega$ is a positive word and $1 \leq i \leq N-2$. We number the bricks of $\tilde\alpha$ in column~$i$ from $1$ to $n$, as before. The braid relation shifts the top brick in column~$i$ to the neighbouring column~$i+1$ of $\tilde\beta$. We keep the label~$n$ for this brick. However, we change the label of the brick~$n-1$ in $\alpha$ to~$n-1'$ in $\beta$, see Figure~7.
\begin{figure}[ht]
\scalebox{0.7}{\raisebox{-0pt}{$\vcenter{\hbox{\includegraphics{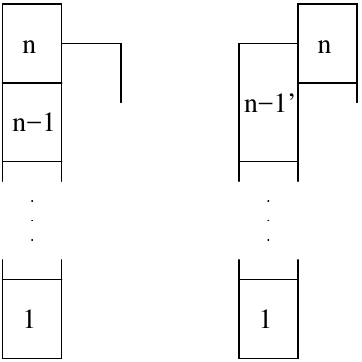}}}$}}
\caption{The braid relation connecting $\alpha$ and $\beta$}
\end{figure}

As in the previous section, we define a homomorphism
$$\varphi: \brickgroupa \to \brickgroup$$
by mapping the generator $s_{n-1}$ to $s_n^{-1} s'_{n-1} s_n$, and leaving all the other generators invariant. We need to show that this map is compatible with all the relations in $R_{\tilde\alpha}$ and $R_{\tilde\beta}$ involving $s_{n-1}$, $s'_{n-1}$, $s_n$. As before, we derive all the relations in $R_{\tilde\beta}$ involving $s'_{n-1}$ or $s_n$ (or both) from relations in $R_{\tilde\alpha}$, by substituting $\varphi^{-1}(s'_{n-1})=s_n s_{n-1} s_n^{-1}$ for $s'_{n-1}$. We start by verifying the braid and commutation relations between $s'_{n-1}$ and the generators $s_1,s_2,\ldots s_{n-2}$, as well as $s_n$.

\begin{enumerate}
\item[(i)] $[s_i,s'_{n-1}]=1$, for all $1 \leq i \leq n-3$.
$$s_i s'_{n-1}=s_i s_n s_{n-1} s_n^{-1}= s_n s_{n-1} s_n^{-1} s_i=s'_{n-1} s_i$$ 

\item[(ii)] $s_{n-2} s'_{n-1} s_{n-2} =s'_{n-1} s_{n-2} s'_{n-1}$.
$$\begin{aligned} s_{n-2} s'_{n-1} s_{n-2} &= s_{n-2} s_n s_{n-1} s_n^{-1} s_{n-2}=s_n s_{n-1} s_{n-2} s_{n-1} s_n^{-1}  \\
&= s_n s_{n-1} s_n^{-1} s_{n-2} s_n s_{n-1} s_n^{-1}=s'_{n-1} s_{n-2} s'_{n-1}
\end{aligned}$$

\item[(iii)] $s_n s'_{n-1} s_n=s'_{n-1} s_n s'_{n-1}$.
$$s_n s'_{n-1} s_n=s_n^2 s_{n-1}=s_n s_{n-1} s_n s_{n-1} s_n^{-1}=s'_{n-1} s_n s'_{n-1}$$
\end{enumerate}

Next, we verify all the relations in $R_{\tilde\beta}$ involving $s'_{n-1}$ and a generator $s_a$ associated with a brick $a$ from another column. As in the previous section, we may assume that the brick~$a$ is located in one of the neighbouring columns~$i+1$ or~$i-1$. There is an easy case, when the bricks~$a$ and~$n$ are not linked, which we omit here. The remaining three cases of how the brick~$a$ is linked with the bricks of column~$i$ are shown in Figure~8.
\begin{figure}[ht]
(1) \, \scalebox{0.7}{\raisebox{-0pt}{$\vcenter{\hbox{\includegraphics{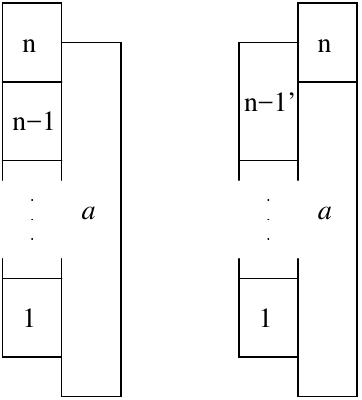}}}$}} \quad
(2) \, \scalebox{0.7}{\raisebox{-0pt}{$\vcenter{\hbox{\includegraphics{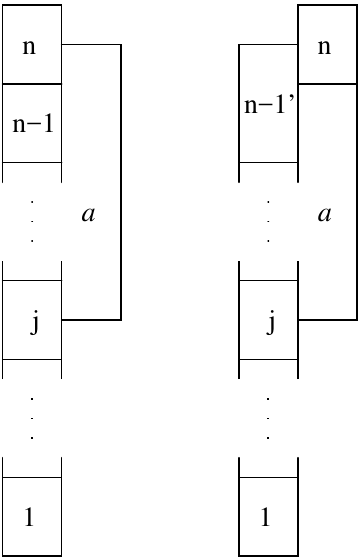}}}$}} \quad
(3) \, \scalebox{0.7}{\raisebox{-0pt}{$\vcenter{\hbox{\includegraphics{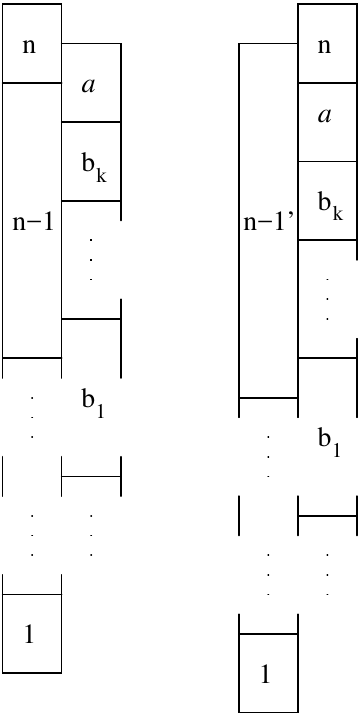}}}$}}
\caption{Three positions for the brick~$a$}
\end{figure}

\begin{enumerate}
\item The brick~$a$ is linked with brick~$n$, and no other brick of column~$i$. We need to derive the braid relation $s_a s'_{n-1} s_a=s'_{n-1} s_a s'_{n-1}$, as well as the cycle relation for the triangle of bricks~$n-1',n,a$.
$$\begin{aligned} s_a s'_{n-1} s_a &= s_a s_n s_{n-1} s_n^{-1} s_a=s_a s_{n-1}^{-1} s_n s_{n-1} s_a \\
&= s_{n-1}^{-1} s_n s_a s_n s_{n-1}=s_{n-1}^{-1} s_n s_{n-1} s_a s_{n-1}^{-1} s_n s_{n-1} \\
&= s'_{n-1} s_a s'_{n-1}.
\end{aligned}$$

$$\begin{aligned} s_n s_a s'_{n-1} s_n &= s_n s_a s_n s_{n-1}=s_a s_n s_a s_{n-1} \\
&= s_a s_n s_{n-1} s_a=s_a s'_{n-1} s_n s_a.
\end{aligned}$$

\item The brick~$a$ is linked with brick~$n$ and one other brick $j \leq n-2$ from column~$i$. In addition to case~(1), we need to derive the cycle relation for the cycle of bricks~$n-1',n-2,\ldots,j+1,j,a$.
For this purpose, we compute
$$\begin{aligned} A &= s_n s_{n-1} \ldots s_j s_a s_n s_{n-1} \ldots s_{j+2} s_{j+1} \\
&= s_n s_{n-1} \ldots s_j s_a s'_{n-1} s_n s_{n-2} \ldots s_{j+2} s_{j+1}, \\
B &= s_{n-1} s_{n-2} \ldots s_j s_a s_n s_{n-1} \ldots s_j \\
&= s_{n-1} s_{n-2} \ldots s_j s_a s'_{n-1} s_n s_{n-2} \ldots s_j, \\
\end{aligned}$$
$$\begin{aligned}
\qquad\qquad C=s_{n-1}^{-1} A s_n^{-1} &= s_{n-1}^{-1} s_n s_{n-1} \ldots s_j s_a s_n s_{n-1} s_n^{-1} s_{n-2} \ldots s_{j+2} s_{j+1} \\
&= s'_{n-1} \ldots s_j s_a s'_{n-1} s_{n-2} \ldots s_{j+2} s_{j+1}, \\
D=s_{n-1}^{-1} B s_n^{-1} &= s_{n-2} \ldots s_j s_a s_n s_{n-1} s_n^{-1} s_{n-2} \ldots s_j \\
&= s_{n-2} \ldots s_j s_a s'_{n-1} s_{n-2} \ldots s_j.
\end{aligned}$$
Thus the cycle relation $C=D$ for the bricks~$n-1',n-2,\ldots,j+1,j,a$ is equivalent to the cycle relation $A=B$ for the bricks~$n,n-1,\ldots,j+1,j,a$.

\item The brick~$a$ is linked with brick~$n$ and~$n-1$. In this case, there is possibly a brick~$b_1$ below the brick~$a$, which is also linked with the brick~$n-1$. We number the bricks between~$b_1$ and~$a$ from the bottom to the top: $b_2,b_3,\ldots,b_k$. We need to derive the commutation relations $[s_a,s'_{n-1}]=1$ and $[s_{b_j},s'_{n-1}]=1$ (which is trivial), as well as the cycle relation for the cycle of bricks~$n-1',n,a,b_k,\ldots,b_1$. For the commutation relation $[s_a,s'_{n-1}]=1$, we compute
$$s_a s'_{n-1} s_n s_a=s_a s_n s_{n-1} s_a,$$
$$s'_{n-1} s_a s_n s_a=s'_{n-1} s_n s_a s_n=s_n s_{n-1} s_a s_n.$$
Thus $[s_a,s'_{n-1}]=1$ follows from the cycle relation for the triangle of bricks~$n,n-1,a$.

Finally, we show that the cycle relations for the two cycles of bricks
$n-1,a,b_k,\ldots,b_1$ and $n-1',n,a,b_k,\ldots,b_1$ are equivalent.
Multiplying the relation from the first cycle of bricks from the left by $s_a s_n$, we obtain the following four equivalent lines:
$$s_a s_n s_a s_{b_k} \ldots s_{b_1} s_{n-1} s_a \ldots s_{b_2}=
s_a s_n s_{b_k} \ldots s_{b_1} s_{n-1} s_a \ldots s_{b_1},$$
$$s_n s_a s_n s_{b_k} \ldots s_{b_1} s_{n-1} s_a \ldots s_{b_2}=
s_a s_n s_{b_k} \ldots s_{b_1} s_{n-1} s_a \ldots s_{b_1},$$
$$s_n s_a s_{b_k} \ldots s_{b_1} s_n s_{n-1} s_a \ldots s_{b_2}=
s_a s_{b_k} \ldots s_{b_1} s_n s_{n-1} s_a \ldots s_{b_1},$$
$$s_n s_a s_{b_k} \ldots s_{b_1} s'_{n-1} s_n s_a \ldots s_{b_2}=
s_a s_{b_k} \ldots s_{b_1} s'_{n-1} s_n s_a \ldots s_{b_1}.$$
The last equation is the cycle relation for the second cycle of bricks.
\end{enumerate} 

\section{Conclusion}

To complete the proof of Theorem~1, we have to show that in the presence of a positive half twist
positive words representing conjugate braids can be transformed into each other by positive braid relations and elementary conjugations only.

Recall first the foundational article of Garside \cite{Gar}. He defines a left normal form for every braid $\beta$,
i.e.\ a uniquely defined $k_\beta\in \mathbb Z$ and a uniquely defined positive braid $\beta'$ such that 
$\beta = \Delta^{k_\beta}\beta'$. The \emph{summit power} is the maximal value $\mathbf k_\beta$ the exponent
takes on the conjugacy class of $\beta$, which is shown to exist.

Note in particular, a braid $\beta$ is positive if $k_\beta$ is non-negative and
its braid isotopy class contains a positive half twist if $\mathbf k_\beta$ positive.

The conjugacy problem is solved, since the \emph{summit set} of a braid $\beta$
\[
S_\beta \quad = \quad \{\: \alpha\sim\beta \:|\: k_\alpha=\mathbf k_\beta, \:\alpha\mbox{ in normal form} \:\}
\]
is a finite set, can be determined by an algorithm in finite time from $\beta$, and is a complete invariant
of the conjugacy class of $\beta$.

Garside proves moreover that positive words representing braids in $S_\beta$ are related through positive words such that
for each subsequent pair $\tilde\alpha,\tilde\alpha'$ of positive words there exists a (possibly trivial) 
\emph{permutation braid} $\gamma$, i.e.\ a positive left divisor of $\Delta$ with a positive representative $\tilde\gamma$ such that
$\tilde\alpha\tilde\gamma$ and $\tilde\gamma\tilde\alpha'$ are related through positive words by braid relations only.
We say in this case, that the two words are related by a \emph{permutation braid conjugation}.

In a similar direction Elrifai and Morton \cite{EM} prove that any positive word $\tilde\beta$ can be transformed into a word
representing an element of $S_\beta$ by \emph{cycling}, which for positive words can be realized by a sequence
through positive words using positive braid relations and elementary conjugations only.

Hence finally the claim of the theorem follows from the following observation of Orevkov which came to our attention by work of Fiedler:

\begin{quote}
Two words representing braids in $S_\beta$ with positive summit power are related through positive words using 
positive braid relations and elementary conjugations only.
\end{quote}

In fact, it suffices to obtain the claim for two such words related by a permutation braid relation.
Let $\tilde\gamma$ be a representative of the permutation word involved and $\tilde\gamma'$ such that $\tilde\gamma\tilde\gamma'$ transforms
into $\Delta$ by positive braid relations only.

By positivity of the summit power the first word has the form $\Delta^\mathbf k\tilde\alpha$ which transforms
into $\tilde\gamma\tilde\gamma'\Delta^{\mathbf k-1}\tilde\alpha$ by positive braid relations only. Using elementary conjugations
we can further transform to $\tilde\gamma'\Delta^{\mathbf k-1}\tilde\alpha\tilde\gamma$ through positive braids. By
assumption the normal of $\tilde\gamma'\Delta^{\mathbf k-1}\tilde\alpha\tilde\gamma$ is the second braid. So by the algorithm
of Garside the two words are related by positive braid relations only.

\medskip
\noindent
Mathematisches Institut, Sidlerstr.~5, CH-3012 Bern, Switzerland

\smallskip
\noindent
\texttt{sebastian.baader@math.unibe.ch}

\medskip
\noindent
Mathematisches Institut, Universit\"atsstr.~30, D-95447 Bayreuth, Germany

\smallskip
\noindent
\texttt{michael.loenne@uni-bayreuth.de}


\begin{thebibliography}{99} 

    
\bibitem{Ba}
    S.~Baader: \emph{Positive braids of maximal signature}, Enseign. Math.~\textbf{59} (2013), no.~3-4, 351--358.
    
\bibitem{BLL}
    S.~Baader, L.~Lewark, L.~Liechti: \emph{Checkerboard graph monodromies}, Enseign. Math.~\textbf{64} (2018), no.~1-2, 65--88.
    
\bibitem{BH}
    J.~S.~Birman, H.~M.~Hilden: \emph{On isotopies of homeomorphisms of Riemann surfaces} Ann. of Math. (2)~\textbf{97} (1973), 424--439.

\bibitem{BS}
E.~Brieskorn, K.~Saito: \emph{Artin-Gruppen und Coxeter-Gruppen}, 
Invent.\ Math.~\textbf{17} (1972), 245--271.

\bibitem{Deh}
    P.~Dehornoy: \emph{On the zeroes of the Alexander polynomial of a Lorenz knot}, Annales de l'Institut Fourier~\textbf{65}, no.~2 (2015), p.~509--548. 


\bibitem{EM}
      E.~A.~Elrifai, H.~R.~Morton: \emph{Algorithms for Positive Braids}, Quart.~J. Math.\ Oxford \textbf{45} (1994), 479--497.

\bibitem{EH}
    J.~Etnyre, J.~Van~Horn-Morris: \emph{Fibered transverse knots and the Bennequin bound}, Int. Math. Res. Not. IMRN 2011, no.~7, 1483--1509.
    
\bibitem{Ga1} A.M.~Gabrielov:
{\sl Dynkin diagrams for unimodular singularities},
Funct.\ Anal.\ Appl.\ 8 (1974), 192--196

\bibitem{Gar} 
      F.~A.~Garside: \emph{The braid group and other groups}, Quart.~J. Math.\ Oxford \textbf{20} (1969), 235--254.
       
\bibitem{GL}
    C.~McA.~Gordon, J.~Luecke: \emph{Knots are determined by their complements}, J.~Amer. Math. Soc.~\textbf{2} (1989), no.~2, 371--415. 

\bibitem{GM}    
    J.~Grant, R.~J.~Marsh: \emph{Braid groups and quiver mutation}, Pacific~J.~Math.~\textbf{290} (2017), 77--116.

\bibitem{HL}
	A.~Hefez, F.~Lazzeri: \emph{The intersection matrix of Brieskorn singularities}, Invent.\ Math.\ 25 (1974), 143--157.

\bibitem{Lab}
Ch.~Labru\`ere: \emph{Generalized Braid Groups and Mapping Class Groups},
J. Knot Theory Ramifications~\textbf{6} (1997), 715--726.


\bibitem{L}
    M.~L\"onne: \emph{Fundamental group of discriminant complements of Brieskorn-Pham polynomials}, C. R. Math. Acad. Sci. Paris~\textbf{345} (2007), no.~2, 93--96.

\bibitem{L_IJM}
    M.~L\"onne: \emph{Braid monodromy of some Brieskorn-Pham singularities},\\ Intern.\ J.\ Math.~\textbf{21} (2010), no.~8, 1047--1070.
    
\bibitem{L_JST}
    M.~L\"onne: \emph{On a discriminant knot group problem of Brieskorn}, J.\ Singul.~\textbf{18} (2018), 455--463.

\bibitem{PV}
    B.~Perron, J.~P.~Vannier: \emph{Groupe de monodromie g\'eom\'etrique des singularit\'es simples}, Math. Ann.~\textbf{306} (1996), no.~2, 231--245.


\bibitem{St}
    J.~R.~Stallings: \emph{Constructions of fibred knots and links}, Algebraic and geometric topology, Proc. Sympos. Pure Math. XXXII (1978), Part~2, 55--60, Amer. Math. Soc., Providence, R.I.

\bibitem{Waj}
B.~Wajnryb: \emph{Artin groups and geometric monodromy}, 
Invent.\ Math.~\textbf{138} (1999), 563--571.





\end{thebibliography}
\end{document}